\documentclass{amsart}
\usepackage{graphicx}
\usepackage{amsmath}
\usepackage{amssymb}
\usepackage[T1]{fontenc}
\usepackage[english]{babel}

\newtheorem{thm}{Theorem}[section]
\newtheorem{cor}[thm]{Corollary}

\newtheorem{prop}[thm]{Proposition}
\newtheorem{example}[thm]{Exemple}
\theoremstyle{definition}
\newtheorem{definition}[thm]{Definition}
\theoremstyle{Remark}
\newtheorem{rem}[thm]{Remark}
\numberwithin{equation}{section}

\newcommand{\K}{\mathbb K}

\newcommand{\A}{\mathcal{A}}
\newcommand{\HH}{\mathcal{H}}
\newcommand{\B}{\mathcal{B}}

\begin{document}

\title[Kaplansky's Constructions Type and Classification of Weak bialgebras ...]{Kaplansky's
Construction Type and Classification \\ of Weak bialgebras and Weak Hopf algebras}%

\author{Z. Chebel and A. Makhlouf}%
\address{Skikda University, Laboratory of  physics,  surface and  interface, Algeria}
\email{zoheir\_chebel1@yahoo.fr}%
\address{ Haute Alsace University, Laboratoire de Math\'{e}matiques, Informatique et
Applications,  6 rue des Fr\`{e}res Lumi\`{e}re 68093 Mulhouse France}
\email{Abdenacer.Makhlouf@uha.fr}


\keywords{Weak bialgebra, Weak Hopf algebra, construction,
classification, automorphisms group}
\subjclass[2000]{ 16W30.}

\begin{abstract}
In this paper, we study  weak bialgebras and  weak Hopf algebras. These algebras form a class wider
than  bialgebras respectively Hopf algebras.  The main results of this paper are Kaplansky's constructions type which
lead to weak bialgebras or weak Hopf algebras starting from a regular  algebra or a bialgebra. Also we provide a classification
of  $2$-dimensional and $3$-dimensional weak bialgebras and weak Hopf algebras. We  determine then
the stabilizer group and the representative of these classes, the action being that of the linear group.
\end{abstract}

\maketitle


\section*{Introduction}
Motivated by quantum symmetry and field algebras, the weak coproduct was introduced first by Mack G. and Schomerus V. in \cite{Mack-Schomerus,Schomerus}. The weak bialgebras were introduced by B\"{o}hm G., Nill F. and Szlach\'{a}nyi K., in \cite{BNS} with a motivation from operator algebras and quantum field theory.  The
weak Hopf algebras, called also quantum groupoids, appeared also in dynamic deformation theory of
quantum groups \cite{Etingof}.    Weak bialgebras and weak Hopf algebras were developed from algebraic point of view and  have been considered
by several authors in various settings (see \cite{espagn1,espagn2,Caenepeel,Caenepeel2,Caenepeel3,Kadison,Nikshych,vainerman1,Schauenburg,
vallin,Vecse,Wisbauer}). The weak bialgebras (resp. weak Hopf algebras)  constitute a class
wider than
 bialgebras (resp. Hopf algebras) where we do not require the conservation of the unit by the multiplication and the multiplicativity of the counit.

The aim of this work is to provide some constructions of finite-dimensional weak bialgebras starting from any algebras.
These constructions are inspired from the Kaplansky's construction of bialgebras (see \cite{Ka73}[Theorem 1.3]). Also, we show that the set of weak bialgebras (resp. weak Hopf algebras) forms an
algebraic variety fibred by
 a linear group action. The orbits under this action
correspond to isomorphic classes.  Therefore, we determine up to isomorphisms all the
$n$-dimensional weak bialgebras and weak Hopf algebras  with $n \leqslant 3$, for which we compute
the corresponding stabilizer subgroups.

In the first Section of the paper we summarize the definitions and the main properties of weak
bialgebras and weak Hopf algebras. Section 2 is dedicated to Kaplansky's construction type, we provide several constructions of weak
bialgebras (resp. weak Hopf algebras) starting from any algebra or bialgebra (resp. Hopf algebra). In the last Section, we establish a classification up to isomorphism of $n$-dimensional weak
bialgebras and weak Hopf algebras  for $n \leqslant 3$. Then we compute their stabilizer groups.

\section{Generalities}
Throughout this paper $\mathbb{K}$ is an algebraically closed field of characteristic 0. In this Section,
we review briefly the algebraic  theory of weak bialgebras and weak Hopf algebras, see \cite{Hopf-Roumains,Kas95,Ma05} for bialgebras and Hopf algebras theory. Let $V$ be a finite-dimensional $\mathbb{K}$-vector
space. In the sequel we use Sweedler's notation for a comultiplication $\Delta$,  that is
$\Delta(x)=\sum_{(1)(2)}{ x_{(1)}\otimes x_{(2)}}$, $ \forall  x \in V$. The summation sign is omitted when there is no ambiguity. 

A \textit{bialgebra}  is a $\K$-vector space $V$ equipped with an algebra structure given by a
multiplication $m$ and a unit  $\eta$ and a
 coalgebra structure given by a comultiplication
$\Delta$ and a counit $\varepsilon$, such that there is a compatibility condition between these two
structures expressed by the fact that $\Delta$ and $\varepsilon$ are algebras morphisms, that is for
$x,y\in V$
$$\Delta(m(x\otimes y))=\Delta(x)\bullet \Delta(y)\quad \text{and}
\quad \varepsilon(m(x\otimes y))=\varepsilon(x)\varepsilon (y).$$ The multiplication $\bullet $ on $V\otimes
V$  is the usual multiplication on tensor product, $$ \left( x\otimes y\right) \bullet \left( x^{\prime
}\otimes y^{\prime }\right) =m \left(x\otimes x^{\prime } \right) \otimes m \left( y\otimes y^{\prime
}\right).
$$
The unit $\eta$ is completely determined by $\eta (1)$, which we denote by $ 1 $. It is  assumed also that the unit $1$ is conserved by the comultiplication, that is $\Delta (1)=1\otimes 1$.
A bialgebra is said to be a \textit{Hopf algebra} if the identity map on $V$ has an inverse for the convolution
product defined by
\begin{equation}\label{conv} f\star g := m
\circ (f\otimes g)\circ \Delta .
\end{equation}
The unit for the convolution product being $\eta\circ\varepsilon$.
 For simplicity, the  multiplication
$m$ is denoted by a dot  when there is  no confusion. In the following, we recall the definition of weak bialgebra.
\begin{definition}\label{DefinitionWB}
 A weak bialgebra is a quintuple $ \mathcal{B}=(V, m, \eta , \Delta ,\varepsilon), $ where
$ m : V\otimes V \rightarrow V $ (multiplication), $\eta :\mathbb{K} \rightarrow V $ (unit), $ \Delta:V \rightarrow  V
\otimes V $ (comultiplication) and $\varepsilon : V \rightarrow \mathbb{K} $ (counit) are linear maps satisfying :

\begin{enumerate}
\item the triple $(V, m, \eta)$ is a unital associative algebra, that is
\begin{equation}\label{mult1}
 m \left( m
\left( x\otimes y\right)\otimes z\right) -m \left( x\otimes m \left( y\otimes z\right) \right) =0 \quad  \forall
x, y, z\in V,
\end{equation}
\begin{equation}\label{mult2}
   m \left(
x\otimes 1\right) =m \left( 1\otimes x\right) =x  \quad \forall x\in V,
\end{equation}
\item the triple $(V, \Delta, \varepsilon)$ is a coalgebra, that is
\begin{equation}\label{comult1}
(\Delta \otimes id) \Delta (x)=(id\otimes\Delta)\Delta (x) \quad \forall x\in V,
\end{equation}
\begin{equation}\label{comult2}
(\varepsilon \otimes id) \Delta (x)=(id \otimes \varepsilon) \Delta(x)=id(x) \quad \forall x\in V,
\end{equation}
\item the compatibility condition is expressed by the the following three identities:
 \begin{equation}\label{compatibility}
     \Delta (m(x\otimes y)) = \sum_{(1)(2)}{ m(x_{(1)}\otimes y_{(1)}) \otimes m(x_{(2)}\otimes
     y_{(2)})}
     \quad \forall x, y \in V,
\end{equation}
\begin{equation}\label{deltaUnit} (\Delta \otimes id)\Delta(1)=(\Delta(1)\otimes
    1)\bullet(1\otimes \Delta(1) )=(1\otimes \Delta(1) )\bullet(\Delta(1)\otimes 1),
\end{equation}
\begin{equation}\label{deltaCoUnit} \varepsilon (m (m (x \otimes y) \otimes z))=\varepsilon
    (m (x \otimes y_ {(1)} ) )\ \varepsilon(m (y_{(2)} \otimes z)) \quad \forall x, y,z \in V.
\end{equation}
\end{enumerate}
\end{definition}

\begin{rem}The condition \eqref{compatibility} means that
$\Delta$ is an algebras homomorphism. But condition \eqref{deltaUnit} shows that $\Delta$ does not
necessarily conserve the unit $1$. If $\Delta (1)=1\otimes 1$ then the condition \eqref{deltaUnit} is
satisfied.

Identity \eqref{deltaCoUnit} is a weak version of the fact that $\varepsilon$ is an algebras morphism  in
the bialgebra case. Indeed, if $\varepsilon$ is an algebras homomorphism then
\begin{eqnarray*}\varepsilon
    (m (x \otimes y_ {(1)} ) )\ \varepsilon(m (y_{(2)} \otimes z))& =&\varepsilon
   ( x) \varepsilon( y_ {(1)} )  \varepsilon(y_{(2)}) \varepsilon( z)\\ \ & =&\varepsilon
   ( x) \varepsilon( y_ {(1)}   \varepsilon(y_{(2)})) \varepsilon( z) =\varepsilon
   ( x) \varepsilon( y) \varepsilon( z)\\ \ & =&
     \varepsilon (m (m (x \otimes y) \otimes z))
\end{eqnarray*}

 When  $\Delta (1)=1\otimes 1$ then one can derive that the counit is  an algebras homomorphism,
 indeed
 $$\varepsilon (m(x \otimes y))=\varepsilon (m (m (x \otimes 1) \otimes
y))=\varepsilon (m (x \otimes 1)) \varepsilon (m (1 \otimes y)) =\varepsilon (x) \varepsilon (y).$$ Then a
bialgebra is always a weak bialgebra.

One may consider a definition of weak bialgebras where counits $\varepsilon$ are algebras homomorphisms.
In this paper we consider  weak bialgebras where  counits $\varepsilon$   are not necessarily
algebras homomorphisms   but satisfy  the identity \eqref{deltaCoUnit}.
\end{rem}

\begin{definition}
A weak Hopf algebra is a sextuple  $\mathcal{H}=(V, m, \eta, \Delta, \varepsilon, S)$, where $(V, m,
\eta, \Delta, \varepsilon )$ is a weak bialgebra and  $S$ is an antipode which is an endomorphism of
$V$ such as :
\\$\forall x\in V$,
\begin{equation}\label{WHopf1}
 m(id\otimes S)\Delta (x)=(\varepsilon \otimes id)\Delta(1)(x\otimes 1),
\end{equation}
\begin{equation}\label{WHopf2}
m(S\otimes id)\Delta (x)=(id\otimes \varepsilon )(1\otimes x)\Delta (1),
\end{equation}
\begin{equation}\label{WHopf3}
m(m\otimes id)(S\otimes id\otimes S)(\Delta \otimes id)\Delta (x)=S(x).
\end{equation}
\end{definition}
\begin{rem}The antipode when it exists is unique and bijective. It is also both algebra and coalgebra
anti-homomorphism.
\end{rem}

\section{Kaplansky's Constructions Type of weak bialgebras}
In this section, we provide  constructions of finite-dimensional weak bialgebras and weak Hopf algebras starting from any algebra or bialgebra. These constructions are inspired by Kaplansky's constructions for bialgebras (see
\cite{Ka73}).

\begin{thm}\label{KplanskyThm}Let $\mathcal{A}$ be any algebra (not necessarily unital) and
$\mathcal{B}$ be the result of adjoining to $\mathcal{A}$ two successive unit elements $e$ and $ 1 $.
On the vector space $\mathcal{B}$ spanned by the vector space $\mathcal{A}$ together with the
generators $\{ 1,\ e\}$, we set
\begin{eqnarray}
\label{K1}\Delta(1)&=&(1-e)\otimes (1-e)+e\otimes e, \\ \label{K2} \Delta(a)&=&a\otimes a \quad \forall a\in
\mathcal{B}\setminus\{1\}, \\ \label{K3}
 \varepsilon(a)&=&1 \quad\forall a\in \mathcal{B}\setminus\{1\}, \\ \label{K4}
  \varepsilon(1)&=&2.
  \end{eqnarray}
   Then $\mathcal{B}$ becomes a weak bialgebra.
\end{thm}
\begin{proof}
The identities \eqref{mult1},\eqref{mult2} are satisfied. In the following we check the remaining  identities in Definition \ref{DefinitionWB}. First, we show that $\Delta$ is coassociative \eqref{comult1}. Let $a\in
\mathcal{B}\setminus\{1\}$, we have
$(\Delta\otimes id)\Delta(a)=a\otimes a\otimes a=(id\otimes \Delta)\Delta(a).$
We have also $(\Delta\otimes id)\Delta(1)=(id\otimes \Delta)\Delta(1) $ since
\begin{eqnarray*}
(\Delta\otimes id)\Delta(1) &=&\Delta(1)\otimes1-\Delta(1)\otimes e-\Delta(e)\otimes 1+2\Delta
(e)\otimes e\\ \ &=&1\otimes1\otimes1-1\otimes e\otimes1-e\otimes1\otimes1+2e\otimes e
\otimes1-1\otimes1\otimes e+1\otimes e\otimes e+
\\
\ & &e\otimes 1\otimes e-2e\otimes e\otimes e-e\otimes e\otimes 1+2e\otimes e\otimes e\\ \
&=&1\otimes 1\otimes 1-1\otimes e\otimes1-e\otimes 1\otimes1+e\otimes e\otimes1-1\otimes
1\otimes e+1\otimes e\otimes e+e\otimes 1\otimes e\\ \ &=&(1-e)\otimes (1-e)\otimes (1-e)+e\otimes
e\otimes e.
\end{eqnarray*}
and on the other hand
\begin{eqnarray*}
(id\otimes \Delta)\Delta(1)&=&1\otimes\Delta(1)-1\otimes \Delta (e)-e\otimes \Delta(1)+2e\otimes
\Delta (e)
\\
\ &=&1\otimes 1\otimes 1-1\otimes e\otimes1-e\otimes 1\otimes1+e\otimes e\otimes1-1\otimes
1\otimes e+ 1\otimes e\otimes e+e\otimes 1\otimes e
\\
\ &=&(1-e)\otimes (1-e)\otimes (1-e)+e\otimes e\otimes e.
\end{eqnarray*}

The identity \eqref{comult2} is also satisfied. Indeed, let $a\in \mathcal{B}\setminus  \{1\}$
\begin{eqnarray*}
 (\varepsilon \otimes
id )\Delta(a)&=&\varepsilon(a)a=a=id (a), \\ (id \otimes \varepsilon)\Delta(a)&=&\varepsilon(a)a=a=id (a)\\ (\varepsilon\otimes id )\Delta(1)&=&\varepsilon(1)1-\varepsilon(1)e-\varepsilon (e) 1+2\varepsilon (e)
e=id (1), \\ (id\otimes \varepsilon)\Delta(1)&=&\varepsilon(1)1-\varepsilon(e)1-\varepsilon (1)
e+2\varepsilon (e) e=id (1).
\end{eqnarray*}
The comultiplication $\Delta$ is  in fact an algebras homomorphism \eqref{compatibility}. Indeed, let $
a_{1},a_{2}\in \mathcal{B}\setminus  \{1\},$ since $a_{1}\cdot a_{2}\in \mathcal{A}$ we have $$\Delta(a_{1})\bullet\Delta(a_{2})=(a_{1}\otimes
a_{1})\bullet(a_{2}\otimes a_{2}) =a_{1}\cdot a_{2}\otimes a_{1} \cdot a_{2}=\Delta (a_{1} \cdot
a_{2}).$$
Also, for $a\in \mathcal{B}\setminus  \{1\}$ we have
\begin{eqnarray*}
\Delta(a)\bullet\Delta(1)&=& (a\otimes a)\bullet(1\otimes1-1\otimes e-e\otimes 1+2\times1\otimes 1)\\ \ & = &
a\otimes a-a\otimes a-a\otimes a+2a\otimes a=a\otimes a=\Delta(a).\end{eqnarray*}

We check the compatibility condition (\ref{deltaUnit}). Indeed, since $e$ and $1-e$ are orthogonal idempotent elements, that is $e\cdot e=e$, $(1-e)\cdot(1-e)=1-e$ and $(1-e)\cdot e=e\cdot (1-e)=0$, we have
\begin{equation*}
(1\otimes \Delta (1))\bullet(\Delta(1)\otimes 1)=(1-e)\otimes (1-e)\otimes (1-e)+e\otimes e\otimes e=(\Delta\otimes id)\Delta(1).
\end{equation*}
And similarly
\begin{equation*}
 (\Delta(1)\otimes 1)\bullet(1\otimes \Delta (1))=(1-e)\otimes (1-e)\otimes (1-e)+e\otimes e\otimes e=(\Delta\otimes id)\Delta(1).
\end{equation*}

Finally, we check that the identity \eqref{deltaCoUnit} is satisfied for any element in $\mathcal{B}.$ Indeed, let $a_{1},a_{2},a_3 \in \mathcal{B}$ with $a_2\in \mathcal{B}\setminus\{1\}$. Since $a_1\cdot a_2 \in \mathcal{B}\setminus\{1\}$ (resp. $a_2\cdot a_3 \in \mathcal{B}\setminus\{1\}$), then $(a_1\cdot a_2)\cdot a_3 \in \mathcal{B}\setminus\{1\}$ (resp. $a_1\cdot (a_2\cdot a_3) \in \mathcal{B}\setminus\{1\}$). Therefore
 $\varepsilon(a_{1} \cdot a_2 \cdot
a_{3})=1$  and  $\varepsilon(a_{1} \cdot a_2) \varepsilon(a_2 \cdot
a_{3})=1.$

Assume now that $a_2=1$, then the left hand side becomes $\varepsilon(a_{1} \cdot 1 \cdot
a_{3})= \varepsilon(a_{1}  \cdot
a_{3})$ and the right hand side writes
$\varepsilon(a_{1} \cdot 1_{(1)})\varepsilon(  1_{(2)}\cdot
a_{3})=\varepsilon(a_{1} \cdot (1-e))\varepsilon(  (1-e)\cdot
a_{3})+\varepsilon(a_{1} \cdot e)\varepsilon(  e\cdot
a_{3}).$
We consider the following particular cases:
\begin{enumerate}
\item $a_1=1$ and $a_3=1$

$\varepsilon(1 \cdot 1_{(1)})\varepsilon(  1_{(2)}\cdot
1)=\varepsilon(1 -e)\varepsilon(1-e)+\varepsilon(e)\varepsilon(e)  =2=\varepsilon(1).$

\item $a_1=1$ and $a_3\neq 1$

$\varepsilon(1 \cdot 1_{(1)})\varepsilon(  1_{(2)}\cdot
a_{3})=\varepsilon( 1-e)\varepsilon(  (1-e)\cdot
a_{3})+\varepsilon( e)\varepsilon(  e\cdot
a_{3})=\varepsilon(a_3)=1.$

\item $a_1\neq 1$ and $a_3=1$

$\varepsilon(a_1 \cdot 1_{(1)})\varepsilon(  1_{(2)}\cdot
1)=\varepsilon(a_{1} \cdot (1-e))\varepsilon(  1-e
)+\varepsilon(a_{1} \cdot e)\varepsilon(  e\cdot
1)=\varepsilon(a_1)=1.$

\item $a_1\neq 1$ and $a_3\neq 1$

$\varepsilon(a_1 \cdot 1_{(1)})\varepsilon(  1_{(2)}\cdot
a_{3})=\varepsilon(a_{1} \cdot (1-e))\varepsilon(  (1-e)\cdot
a_{3})+\varepsilon(a_{1} \cdot e)\varepsilon(  e\cdot
a_{3})=\varepsilon(a_1)\varepsilon(a_3)=1,$ which is equal to $\varepsilon(a_1\cdot a_3)$ because $a_1\cdot a_3\in \mathcal{B}\setminus\{1\}.$
\end{enumerate}
This ends the proof that $\mathcal{B}$ is endowed with a weak bialgebra structure.
\end{proof}

\begin{rem}The weak bialgebra obtained above is not a regular bialgebra since for $a\in \mathcal{B}\setminus  \{1\}$, we have
$
\varepsilon (a\cdot 1)=1
$ and $
\varepsilon (a)\varepsilon( 1)=2.
$
\end{rem}

\begin{cor}
Let $\A$ be an associative algebra with a unit $ 1 $. If $\A\setminus\{1\}$ is a subalgebra and there
 exists an element $e$ in $\A\setminus\{1\}$ such that $e\cdot e=e$ and
  $e \cdot a = a \cdot e=a$ for all $a\in \A\setminus\{1\}$.
  Then there exists a weak bialgebra structure on $\A$ given by:
\begin{eqnarray*}
\Delta(1)&=&(1-e)\otimes (1-e)+e\otimes e,\\ \Delta(a)&=&a\otimes a, \ \ \forall a\in \A\setminus\{1\},
\\ \varepsilon(1)&=&2, \\ \varepsilon(a)&=&1,\ \ \forall a\in \A\setminus\{1\}.
\end{eqnarray*}

\end{cor}

\begin{rem}The weak bialgebra $\mathcal{B}$ of theorem \ref{KplanskyThm} is not always a weak Hopf algebra with the  antipode $S=id$.
 The identities \eqref{WHopf1},\eqref{WHopf2},\eqref{WHopf3}  are fulfilled for $x=1$.  But for $x\in \mathcal{B}\setminus \{1\}$, the identities \eqref{WHopf1},\eqref{WHopf2} lead to the condition $x\cdot x=e$ while the identity \eqref{WHopf2} leads to $(x\cdot x)\cdot x=x.$
\end{rem}

\begin{example}\label{example1}
Let $\A$ be a 2-dimensional  associative algebra  with unit $1$. We assume that
 $e$ is an idempotent element ($e\cdot e=e$) different from $1$ in $\A$. Then there exists a weak  bialgebra structure on $\A$ given by
\begin{eqnarray*}
\Delta(1)&=&(1-e)\otimes (1-e)+e\otimes e,\\\Delta(e)&=&e\otimes e, \\ \varepsilon(1)&=&2, \quad
\varepsilon(e)=1.
\end{eqnarray*}

Moreover, with  $S=id$,   the bialgebra $\A$ becomes a weak Hopf algebra.
\end{example}

\begin{example}
Let $\A=\K \times \cdots \K$ be an $n$-dimensional  unital associative algebra. Assume that
$\mathfrak{b}=\{e_{i}\}_{1\leq i \leq n}$ is a basis of   $\A$ such that $ e_1=1 $ and  $\{e_{i}\}_{2\leq
i \leq n}$ are orthogonal idempotent elements.  Then there
exist weak bialgebra structures  on $\A$ given by setting, for  a fixed integer $k\in \{2, \cdots n\}, $
\begin{eqnarray*}
 \Delta(1)&=&(1-e_{k})\otimes (1-e_{k})+e_{k}\otimes e_{k},\\
\Delta(e_{i})&=&e_{i}\otimes e_{i}, \ \ \ i\in \{2, \cdots n\}, \\ \varepsilon(1)&=&2, \\ \varepsilon(e_{i})&=&1, \ \ \
i\in \{2, \cdots n\}.
\end{eqnarray*}

\end{example}

In the following, we have the following more general result.

\begin{thm} \label{Thm2.7} Let $\A$ be a finite-dimensional  unital associative  algebra with unit
$e_{1}=1$. Let  $\mathfrak{b}=\mathfrak{b}_{1}\cup \mathfrak{b}_{2}$ be a basis of $\A$  with
$\mathfrak{b}_{1}=\{e_{i}\}_{i=1,\cdots p} $.  Assume that $span(\mathfrak{b}_{2}) $ is  a subalgebra of $\A$ and
\begin{eqnarray*}
&& e_{i} \cdot e_{j}=e_{max(i,j)}, \ \
 \forall i, j=1,\cdots, p,\\
 && e_{i}\cdot f=f\cdot e_{i}=f,\ \ \forall f\in \mathfrak{b}_{2}.
 \end{eqnarray*}
Then the comultiplication $\Delta$ and the counit $\varepsilon$ defined by
\begin{eqnarray*}
&& \Delta(e_{p})=e_{p}\otimes e_{p} \\ && \Delta(e_{i})=(e_{i}-e_{i+1})\otimes
(e_{i}-e_{i+1})+\Delta(e_{i+1}), \ \ \forall i, \ i=1,\cdots, p-1, \\ &&\Delta(f)=f\otimes f , \ \ \forall f\in
\mathfrak{b}_{2}, \\ && \varepsilon(e_{i})=p-i+1, \ \ \forall i,\    i=1,\cdots, p, \\ && \varepsilon(f)=1, \ \
\forall f\in \mathfrak{b}_{2},
\end{eqnarray*}
endow $\A$ with a  weak bialgebra structure.

\begin{proof}

The image  by $\Delta$ of an element  $e_{p-i}\in\mathfrak{b}_{1}$, $i=1,\cdots , p-1$,  can written in the form:\\
$\Delta(e_{i})=(e_{i}-e_{i+1})\otimes(e_{i}-e_{i+1})+
(e_{i+1}-e_{i+2})\otimes(e_{i+1}-e_{i+2})+ \cdots
+(e_{p-1}-e_{p})\otimes(e_{p-1}-e_{p})+e_{p}\otimes e_{p}.$\\

It follows for $i=1,\cdots, p-1$ that $\Delta(e_{i}-e_{i+1})=(e_{i}-e_{i+1})\otimes(e_{i}-e_{i+1})$ and $\varepsilon(e_{i}-e_{i+1})=1.$\\

Let us show that $\Delta $ is coassociative. For $i=1,\cdots, p-1$, we have

$(id\otimes\Delta)\Delta(e_{i})= (id\otimes\Delta)[(e_{i}-e_{i+1})\otimes
(e_{i}-e_{i+1}) \cdots +(e_{p-1}-e_{p})\otimes
(e_{p-1}-e_{p})+e_{p}\otimes e_{p}]=(e_{i}-e_{i+1})\otimes
(\Delta(e_{i})-\Delta(e_{i+1}))+\cdots +(e_{p-1}-e_{p})\otimes
(\Delta(e_{p-1})-\Delta(e_{p}))+e_{p}\otimes\Delta( e_{p})=(e_{i}-e_{i+1})\otimes
(e_{i}-e_{i+1})\otimes(e_{i}-e_{i+1})+ \cdots +(e_{p-1}-e_{p})\otimes
(e_{p-1}-e_{p})\otimes(e_{p-1}-e_{p})+e_{p}\otimes e_{p}\otimes e_{p},$\\

$(\Delta\otimes id)\Delta(e_{i})= (\Delta\otimes id)[(e_{i}-e_{i+1})\otimes
(e_{i}-e_{i+1})+ \cdots +(e_{p-1}-e_{p})\otimes
(e_{p-1}-e_{p})+e_{p}\otimes e_{p}]=(\Delta(e_{i})-\Delta(e_{i+1}))\otimes
(e_{i}-e_{i+1})+ \cdots +(\Delta(e_{p-1})-\Delta(e_{p}))\otimes (e_{p-1}-e_{p})+\Delta(e_{p})\otimes
e_{p}=(e_{i}-e_{i+1})\otimes
(e_{i}-e_{i+1})\otimes(e_{i}-e_{i+1})+ \cdots +(e_{p-1}-e_{p})\otimes
(e_{p-1}-e_{p})\otimes(e_{p-1}-e_{p})+e_{p}\otimes e_{p}\otimes e_{p} .$\\

Then, $(id\otimes\Delta)\Delta(e_{i})=(\Delta\otimes id)\Delta(e_{i}).$ Obviously, one gets the coassociativity for $e_p$ and any $f\in \mathfrak{b}_{2}.$\\

We show that $\varepsilon$ is a counit. For $i=1,\cdots, p-1$, we have

$(id\otimes\varepsilon)\Delta(e_{i})= (id\otimes\varepsilon)[(e_{i}-e_{i+1})\otimes
(e_{i}-e_{i+1})+ \cdots +(e_{p-1}-e_{p})\otimes
(e_{p-1}-e_{p})+e_{p}\otimes e_{p}]=(e_{i}-e_{i+1})(\varepsilon (e_{i})-\varepsilon
(e_{i+1}))+ \cdots +(e_{p-1}-e_{p})(\varepsilon (e_{p-1})-\varepsilon(e_{p}))+e_{p}\varepsilon(e_{p})=
e_{i}-e_{i+1}+e_{i+1}-e_{i+2}+e_{i+2}+ \cdots -e_{p-1}+e_{p-1} -e_{p}+e_{p}=id(e_{i}),$\\

$(\varepsilon\otimes id)\Delta(e_{i})= (\varepsilon\otimes id)[(e_{i}-e_{i+1})\otimes
(e_{i}-e_{i+1})+ \cdots +(e_{p-1}-e_{p})\otimes
(e_{p-1}-e_{p})+e_{p}\otimes e_{p}]=(\varepsilon(e_{i})-\varepsilon (e_{i+1}))
(e_{i}-e_{i+1})+ \cdots +(\varepsilon(e_{p-1})-\varepsilon(e_{p}))
(e_{p-1}-e_{p})+\varepsilon(e_{p})
e_{p}=e_{i}-e_{i+1}+e_{i+1}-e_{i+2}+e_{i+2}+ \cdots -e_{p-1}+e_{p-1} -e_{p}+e_{p}=id(e_{i}).$
\\

Then $(id\otimes\varepsilon)\Delta(e_{i})=(\varepsilon\otimes id)\Delta(e_{i})=id(e_{i}).$ The coassociativity is obviously satisfied for the grouplike elements. \\

The comultiplication $\Delta $ is compatible with the multiplication. Indeed, let $e_{i},e_{k}\in \mathfrak{b}_{1}$, $i,k=1 \cdots p-1,$ and $f\in \mathfrak{b}_{2},$\\
$\Delta(e_{i})\bullet\Delta(f)=[(e_{i}-e_{i+1})\otimes(e_{i}-e_{i+1})+ \cdots
+(e_{p-1}-e_{p})\otimes(e_{p-1}-e_{p})+e_{p}\otimes e_{p}].(f\otimes f) =e_{p}\cdot f\otimes
e_{p}\cdot f =f\otimes f=\Delta(e_{i}\cdot f). $

Similarly we have $\Delta(f)\bullet\Delta(e_{i})=\Delta(f\cdot e_{i}).$ \\

Assume   $i\geq k$, by a direct  calculation we have \\ $\Delta(e_{i})\bullet\Delta(e_{k})=[(e_{i}-e_{i+1})
\otimes(e_{i}-e_{i+1})+ \cdots
+(e_{p-1}-e_{p})\otimes(e_{p-1}-e_{p})+e_{p}\otimes
e_{p}]\bullet[(e_{k}-e_{k+1})\otimes(e_{k}-e_{k+1})+
  \cdots
+(e_{p-1}-e_{p})\otimes(e_{p-1}-e_{p})+e_{p}\otimes e_{p}]=(e_{i}-e_{i+1}) \otimes(e_{i}-e_{i+1})+\cdots
+(e_{p-1}-e_{p})\otimes(e_{p-1}-e_{p})+e_{p}\otimes e_{p}= \Delta(e_{i})=\Delta(e_{i}\cdot
e_{k}).$

Also  for any $f_1,f_2\in \mathfrak{b}_{2}$, we have $\Delta(f_{1})\bullet\Delta(f_{2})=(f_{1}\otimes
f_{1})\bullet(f_{2}\otimes f_{2}) =f_{1}\cdot f_{2}\otimes f_{1} \cdot f_{2}=\Delta (f_{1} \cdot
f_{2}).$\\

In the following we check the identities \eqref{deltaUnit}. We have

$[\Delta (e_{1})\otimes e_{1}]\cdot[e_{1}\otimes \Delta(e_{1})]= [
(e_{1}-e_{2})\otimes(e_{1}-e_{2})\otimes e_{1}+  \cdots
  + (e_{p-1}-e_{p})\otimes(e_{p-1}-e_{p})\otimes
e_{1}+ e_{p}\otimes e_{p}\otimes e_{1}] \cdot[e_{1}\otimes
(e_{1}-e_{2})\otimes(e_{1}-e_{2})+  \cdots  +e_{1}\otimes
(e_{p-1}-e_{p})\otimes(e_{p-1}-e_{p})+e_{1}\otimes e_{p}\otimes e_{p}]=
(e_{1}-e_{2})\otimes(e_{1}-e_{2})\otimes (e_{1}-e_{2})+ \cdots  +
(e_{p-1}-e_{p})\otimes(e_{p-1}-e_{p})\otimes (e_{p-1}-e_{p})+ e_{p}\otimes e_{p}\otimes e_{p}, $

and

$[e_{1}\otimes
\Delta(e_{1})]\cdot[\Delta (e_{1})\otimes e_{1}]= [e_{1}\otimes
(e_{1}-e_{2})\otimes(e_{1}-e_{2})+ \cdots +e_{1}\otimes (e_{p-1}-e_{p})\otimes(e_{p-1}-e_{p})+
e_{1}\otimes e_{p}\otimes e_{p}].[ (e_{1}-e_{2})\otimes(e_{1}-e_{2})\otimes e_{1}+
 \cdots   +
(e_{p-1}-e_{p})\otimes(e_{p-1}-e_{p})\otimes e_{1}+ e_{p}\otimes e_{p}\otimes e_{1}]=
(e_{1}-e_{2})\otimes(e_{1}-e_{2})\otimes (e_{1}-e_{2})+ \cdots   +
(e_{p-1}-e_{p})\otimes(e_{p-1}-e_{p})\otimes (e_{p-1}-e_{p})+ e_{p}\otimes e_{p}\otimes e_{p}. $

Then $[e_{1}\otimes
\Delta(e_{1})]\bullet[\Delta (e_{1})\otimes e_{1}]=[\Delta (e_{1})\otimes e_{1}]\bullet[e_{1}\otimes \Delta(e_{1})]=(\Delta\otimes id)\Delta(e_{1}).$\\

Now we check the identity \eqref{deltaCoUnit}. We consider first a triple $(e_i,e_j,e_k)$. Assume that $ j\leq k $, in this case  $\varepsilon(e_{i}\cdot e_{j}\cdot
e_{k})=\varepsilon(e_{i}\cdot e_{k}) =p-\max(i,k)+1. $\\

On the other hand we have

$\varepsilon(e_{i}\cdot(e_{j})_{(1)})\varepsilon((e_{j})_{(2)}\cdot e_{k})
=\varepsilon(e_{i}\cdot(e_{j}-e_{j+1}))\varepsilon((e_{j}-e_{j+1})\cdot e_{k}))+
 \cdots +\varepsilon(e_{i} \cdot(e_{p-1}-e_{p}))\varepsilon((e_{p-1}-e_{p})\cdot e_{k})+
\varepsilon(e_{i}\cdot e_{p})\varepsilon((e_{p}\cdot e_{k})=
\varepsilon(e_{i}\cdot(e_{k}-e_{k+1}))\varepsilon(e_{k}-e_{k+1})+
 \cdots +\varepsilon(e_{i}\cdot(e_{p-1}-e_{p}))\varepsilon(e_{p-1}-e_{p})+ \varepsilon(e_{i}\cdot
e_{p})\varepsilon(e_{p})= \varepsilon(e_{i}\cdot(e_{k}-e_{k+1}))+ \varepsilon(e_{i}\cdot
(e_{k+1}-e_{k+2}))+  \cdots +\varepsilon(e_{i}\cdot(e_{p-1}-e_{p}))+ \varepsilon(e_{i}\cdot
e_{p})=\varepsilon(e_{i}\cdot e_{k})- \varepsilon(e_{i}\cdot e_{k+1})+ \varepsilon(e_{i}\cdot
e_{k+1})-\varepsilon(e_{i}\cdot e_{k+2})+ \varepsilon(e_{i}\cdot e_{k+2})+  \cdots -\varepsilon(e_{i}\cdot
e_{p-1})+\varepsilon(e_{i}\cdot e_{p-1}) -\varepsilon(e_{i}\cdot e_{p})+\varepsilon(e_{i}\cdot e_{p})=
\varepsilon(e_{i}\cdot e_{k})=p-\max(i,k)+1.$\\

If $j>k$, then  $\varepsilon(e_{i}\cdot e_{j}\cdot e_{k})=\varepsilon(e_{i}\cdot e_{j}) =p-\max(i,j)+1.
$\\

Also

$\varepsilon(e_{i}\cdot(e_{j})_{(1)})\varepsilon((e_{j})_{(2)}\cdot e_{k})
=\varepsilon(e_{i}\cdot(e_{j}-e_{j+1}))\varepsilon((e_{j}-e_{j+1})\cdot e_{k})+
 \cdots +\varepsilon(e_{i}\cdot(e_{p-1}-e_{p})\varepsilon((e_{p-1}-e_{p})\cdot e_{k})+
\varepsilon(e_{i}\cdot e_{p})\varepsilon(e_{p}\cdot e_{k})=
\varepsilon(e_{i}\cdot(e_{j}-e_{j+1}))\varepsilon(e_{j}-e_{j+1})+
\varepsilon(e_{i}\cdot(e_{j+1}-e_{j+2}))\varepsilon(e_{j+1}-e_{j+2})+
 \cdots +\varepsilon(e_{i}\cdot(e_{p-1}-e_{p}))\varepsilon(e_{p-1}-e_{p})+ \varepsilon(e_{i}\cdot
e_{p})\varepsilon(e_{p})=
\varepsilon(e_{i}\cdot(e_{j}-e_{j+1}))+
 \cdots +\varepsilon(e_{i}\cdot(e_{p-1}-e_{p}))+ \varepsilon(e_{i}\cdot e_{p})=\varepsilon(e_{i}\cdot e_{j})-
\varepsilon(e_{i}\cdot e_{j+1})+ \varepsilon(e_{i}\cdot e_{j+1})+  \cdots -\varepsilon(e_{i}\cdot e_{p-1})+\varepsilon(e_{i}\cdot e_{p-1})
-\varepsilon(e_{i}\cdot e_{p})+\varepsilon(e_{i}\cdot e_{p})= \varepsilon(e_{i}\cdot
e_{j})=p-\max(i,j)+1.$\\

For a triple $(f,e_i,e_j)$, we obtain $\varepsilon(f\cdot e_{i}\cdot e_{j})
=\varepsilon(f) =1$ and on the other hand

$\varepsilon(f\cdot(e_{i})_{(1)})\varepsilon((e_{i})_{(2)}\cdot e_{j})
=\varepsilon(f\cdot(e_{i}-e_{i+1}))\varepsilon((e_{i}-e_{i+1})\cdot e_{j}))+
 \cdots +\varepsilon(f \cdot(e_{p-1}-e_{p}))\varepsilon((e_{p-1}-e_{p})\cdot e_{j})+
\varepsilon(f\cdot e_{p})\varepsilon((e_{p}\cdot e_{j})=
\varepsilon(f\cdot e_{p})\varepsilon(e_{p}\cdot e_j )=
\varepsilon(f)\varepsilon(e_{p} )=1.$\\

For a triple $(a_1,f,a_2)$ where $a_1,a_2\in \A$ we obtain $\varepsilon(a_1\cdot f\cdot a_2)
=\varepsilon(f) =1$ because $a_1\cdot f$ and $f\cdot a_2$ belong to $span (\mathfrak{b}_{2})$, which is in fact an ideal. On the other hand we have $\varepsilon(a_{1}\cdot(f)_{(1)})\varepsilon((f)_{(2)}\cdot a_{2})=\varepsilon(f)\varepsilon( f)=1.$

\end{proof}

\end{thm}

We show in the following that the commutative algebra  $\A=\K \times \cdots \K$ carries a structure of weak Hopf algebra. To this end we write the algebra in a suitable basis.
\begin{prop}\label{propo2.8}Let $\A$ be a unital algebra with unit $e_{2}$ such that on
  on a basis $\{e_i\}_{i=2,\cdots,n}$ of $\A$  the multiplication is given by $m(e_{i}\otimes
e_{j})=e_{max(i,j)},\ \ i, j=2,\cdots ,n$. Let $\B$ be the  result of adjoining a second
unit $e_{1}=1 $ to $\A$.
Set
\begin{eqnarray*}
\Delta(e_{n})&=&e_{n}\otimes e_{n},\\ \Delta(e_{i})&=&(e_{i}-e_{i+1})\otimes
(e_{i}-e_{i+1})+\Delta(e_{i+1}), \quad\forall i,\  i=1 \cdots n-1,\\ \varepsilon(e_{i})&=&n-i+1, \quad \forall i,\
i=1 \cdots n, \\ S&=&id
\end{eqnarray*}
Then $\B$ becomes a weak Hopf algebra.
\end{prop}

\begin{proof} The structure of weak bialgebra follows from the previous theorem. It remains to verify
 the antipode's identities \eqref{WHopf1}\eqref{WHopf2}\eqref{WHopf3}.  We have for $i=1,\cdots,n-1$\\

 $m(id\otimes S)\Delta(e_{i})=m(\Delta(e_{i}))=m((e_{i}-e_{i+1})\otimes(e_{i}-e_{i+1})+
 \cdots +
(e_{n-1}-e_{n})\otimes(e_{n-1}-e_{n})+e_{n}\otimes e_{n})=e_{i}-e_{i+1}+e_{i+1}+ \cdots -e_{n-1}+ e_{n-1}-e_{n}+e_{n}=e_{i},$\\

$(\varepsilon\otimes id)[\Delta(e_{1})\bullet(e_{i}\otimes e_{1})]=(\varepsilon\otimes
id)([(e_{1}-e_{2})\otimes(e_{1}-e_{2})+\cdots
  +(e_{n-1}-e_{n})\otimes(e_{n-1}-e_{n})+e_{n}\otimes
e_{n}]\cdot(e_{i}\otimes e_{1}))=\varepsilon(e_{i}-e_{i+1})(e_{i}-e_{i+1})+
 \cdots + \varepsilon(e_{n-1}-e_{n})(e_{n-1}-e_{n})+
\varepsilon(e_{n})e_{n}=e_{i}-e_{i+1}+e_{i+1}+ \cdots -e_{n-1}+ e_{n-1}-e_{n}+e_{n}=e_{i}, $\\

and similarly  $(id\otimes \varepsilon)[(e_{1}\otimes e_{i})\bullet\Delta(e_{1})]=e_{i}. $\\

Thus \eqref{WHopf1}\eqref{WHopf2} hold.
 \\

The identity \eqref{WHopf3} is also satisfied, we use a previous calculation of $(\Delta\otimes id)(\Delta(e_{i}))$ and  $m((e_{i}-e_{i+1})\otimes(e_{i}-e_{i+1}))=e_{i}-e_{i+1}$, then

 $m(m\otimes id)(S\otimes id \otimes S)(\Delta\otimes id)(\Delta(e_{i}))=m(m\otimes id)(id\otimes id \otimes id)(\Delta\otimes id)(\Delta(e_{i}))=m(m\otimes id)(\Delta\otimes id)(\Delta(e_{i}))=
m(m\otimes id)((e_{i}-e_{i+1})\otimes(e_{i}-e_{i+1})+
 \cdots +
(e_{n-1}-e_{n})\otimes(e_{n-1}-e_{n})+e_{n}\otimes e_{n})=m((e_{i}-e_{i+1})\otimes(e_{i}-e_{i+1})+
 \cdots +
(e_{n-1}-e_{n})\otimes(e_{n-1}-e_{n})+e_{n}\otimes e_{n})=e_{i}=S(e_{i}). $

The proof for $e_n$ leads from easy direct calculations.
\end{proof}

%
%
%
%

\begin{rem}
Similarly we may endow the algebra generated by $n$ orthogonal idempotent
elements  by a structure of weak Hopf algebra. The algebra structure is isomorphic to the algebra structure considered in Proposition \ref{propo2.8}, one may consider the same comultiplication  and counit as in this Proposition.
\end{rem}

Now, we provide constructions of  weak bialgebras starting from any bialgebra. We show that any $n$-dimensional bialgebra could be extended to $(n+1)$-dimensional weak bialgebra.

\begin{thm}\label{ConstBialg}
 Let $\B$ be a bialgebra and  $e_{2}$ its unit. We consider the
  set $\mathcal{B}'$ as a result of  adjoining a unit $e_{1}$ to $\B$
  with respect to the multiplication.

  Assume
 \begin{eqnarray*}
&&\Delta(e_{1})=(e_{1}-e_{2})\otimes (e_{1}-e_{2})+e_{2}\otimes e_{2},
 \\ &&\varepsilon(e_{1})=2 .\end{eqnarray*} Then $\mathcal{B}'$ becomes a weak bialgebra.
\end{thm}

\begin{proof}

The identities \eqref{mult1}-\eqref{deltaUnit} follow from the proof of Theorem \ref{KplanskyThm} and the fact that $\Delta(e_{2})=e_{2}\otimes e_{2}$ and $\varepsilon(e_{2})=1$ since $\B$ is a bialgebra with unit $e_2$. It remains to check
the compatibility of the counit with the comultiplication. The identity \eqref{deltaCoUnit} is satisfied when it deals with 3 elements of $\B$.\\

For a triple $(e_{1} , a, b)$, we have
$$\varepsilon (e_1\cdot a_{(1)})\varepsilon (a_{(2)}\cdot b)=\varepsilon ( a_{(1)})\varepsilon (a_{(2)}\cdot b)=\varepsilon (\varepsilon ( a_{(1)})a_{(2)}\cdot b)=\varepsilon ( a\cdot b)=\varepsilon ( e_1\cdot a\cdot b).
$$
The case of triples $(a, b,e_{1})$ is similar. Let us consider now triples of the form $( a,e_{1} , b)$.
The left hand side of \eqref{deltaCoUnit} becomes $\varepsilon(a \cdot e_{1} \cdot
b)= \varepsilon(a  \cdot b)$ and the right hand side writes
$\varepsilon(a \cdot {e_1}_{(1)})\varepsilon(   {e_1}_{(2)}\cdot
b)=\varepsilon(a \cdot (e_1-e_2))\varepsilon(  (e_1-e_2)\cdot
b)+\varepsilon(a \cdot e_2)\varepsilon(  e_2\cdot
b).$
We consider the following particular cases:
\begin{enumerate}
\item $a=e_1$ and $b=e_1$

$\varepsilon(e_1 \cdot {e_1}_{(1)})\varepsilon(  {e_1}_{(2)}\cdot
e_1)=\varepsilon(e_1 -e_2)\varepsilon(e_1-e_2)+\varepsilon(e_2)\varepsilon(e_2)  =2=\varepsilon(e_1).$

\item $a=1$ and $b\neq e_1$

$\varepsilon(e_1 \cdot {e_1}_{(1)})\varepsilon(  {e_1}_{(2)}\cdot
b)=\varepsilon( e_1-e_2)\varepsilon(  (e_1-e_2)\cdot
b)+\varepsilon( e_2)\varepsilon(  e_2\cdot
b)=\varepsilon(b).$

\item $a\neq e_1$ and $b=e_1$

$\varepsilon(a \cdot {e_1}_{(1)})\varepsilon(  {e_1}_{(2)}\cdot
e_1)=\varepsilon(a \cdot (e_1-e_2))\varepsilon(  e_1-e_2
)+\varepsilon(a \cdot e_2)\varepsilon(  e_2\cdot
e_1)=\varepsilon(a).$

\item $a\neq 1$ and $b\neq 1$

$\varepsilon(a \cdot {e_1}_{(1)})\varepsilon(  {e_1}_{(2)}\cdot
b)=\varepsilon(a \cdot (e_1-e_2))\varepsilon(  (e_1-e_2)\cdot
b)+\varepsilon(a \cdot e_2)\varepsilon(  e_2\cdot
b)=\varepsilon(a)\varepsilon(b)=\varepsilon(a\cdot b)$ because $a, b\in \mathcal{B}$ and $\mathcal{B}$ is a bialgebra.
\end{enumerate}

\end{proof}
The dimension of $\B '$ is $dim \mathcal{B}'=dim \B+1$.
\begin{rem}
The counit of $\mathcal{B}$ is not an algebras homomorphism, indeed $\varepsilon (e_1\cdot e_2)=\varepsilon ( e_2)=1$ while $\varepsilon (e_1)\varepsilon ( e_2)=2.$
\end{rem}

The following Theorem provides a way for extending an $n$-dimensional Hopf algebra to an $(n+1)$-dimensional weak Hopf algebra.
\begin{thm}\label{ConstHopf}
 Let $\HH$ be a Hopf algebra and  $e_{2}$ its unit. We consider the
  set $\mathcal{H}'$ as a result of  adjoining a unit $e_{1}$ to $\HH$
  with respect to the multiplication.

  Assume
 \begin{eqnarray*}
&&\Delta(e_{1})=(e_{1}-e_{2})\otimes (e_{1}-e_{2})+e_{2}\otimes e_{2},
 \\ &&\varepsilon(e_{1})=2 \\&& S(e_1)=e_1.\end{eqnarray*} Then $\mathcal{H}'$ becomes a weak Hopf algebra.
\end{thm}
\begin{proof} The structure of weak bialgebra follows from Theorem \ref{ConstBialg}. The remaining identities are given for $e_1$ by straightforward calculations.
\end{proof}

\begin{example}[\textbf{Sweedler's 5-dimensional weak Hopf algebra}]
Assume that $char(\K)\neq 2$. Let $\HH$ be the Sweedler's 4-dimensional Hopf algebra given by
generators and relations as follows: $\HH$ is generated as a $\K$-algebra by c and x satisfying the
relations:
\begin{equation}\label{SweedlerHopf}
c^{2}=e,\  x^{2}=0,\ x\cdot c=-c\cdot x, \quad \text{where } e \text{ being the unit}.
\end{equation}
Let $\HH'$ be the algebra obtained by  adjoining a new unit $ 1
$ to $\HH$. Then $\HH'$ is 5-dimensional  weak bialgebra defined as a $\mathbb{K}$-algebra, with basis a $\{1, e, x, c, c\cdot x\}$ and relations \eqref{SweedlerHopf} and $e\cdot c=c\cdot  e=c,\ e\cdot x=x\cdot e=x,\ e\cdot e=e.$
The coalgebra structure is defined by :
\begin{eqnarray*}
&&\Delta(1)=(1-e)\otimes (1-e)+e\otimes e, \\ &&\Delta(c)=c\otimes c,  \ \Delta(e)=e\otimes e, \
\Delta(x)=c\otimes x+x\otimes e,\\ &&\varepsilon(1)=2,\ \varepsilon(e)=1, \ \varepsilon(c)=1, \
\varepsilon(x)=0,
\end{eqnarray*}
The antipode is given by:
\begin{equation}
S(1)=1,\  S(e)=e,\  S(c)=c,\  S(x)=-c \cdot x.
\end{equation}
 This weak Hopf algebra is non-commutative and non-cocommutative.

\end{example}
\begin{example}[\textbf{Taft's weak Hopf algebras}]
Let $n\geq 2$ be an integer and $\lambda$ be a primitive $n$-th root of unity. Consider the Taft's algebras
$\mathcal{H}_{n^{2}}(\lambda)$, generalizing Sweedler's Hopf algebra,  defined by the generators $c$ and $x$ and where $e$ being the unit, with the relations:
\begin{equation}
c^{n}=e,\  x^{n}=0,\ x\cdot
c=\lambda\ c\cdot x .
\end{equation}
Let $\mathcal{H}'$ be the algebra obtained by adjoining a new  unit $1$ to
$\mathcal{H}_{n^{2}}(\lambda)$.

We set a coalgebra structure defined by:
\begin{eqnarray*}
&&\Delta(1)=(1-e)\otimes(1-e)+e\otimes e,\\ &&\Delta(e)=e\otimes e, \ \Delta(c)=c\otimes c, \
\Delta(x)=c\otimes x+x\otimes e, \\ && \varepsilon(1)=2, \ \varepsilon(e)=1, \ \varepsilon(c)=1,\
\varepsilon(x)=0.
\end{eqnarray*}
Then $\mathcal{H}'$  becomes an  $(n^{2}+1)$-dimensional weak bialgebra, having a basis $\{1, c^{i}x^{j},0\leq i,j\leq
n-1 \}.$\\ It carries a structure of weak Hopf algebra with an  antipode  defined by:
\begin{equation}
S(1)=1,\ S(e)=e,\  S(c)=c^{-1},\ S(x)=-c^{-1}\cdot x.
\end{equation}

\end{example}

Next two propositions give other constructions of weak bialgebras starting from  bialgebras, the proofs are similar to previous ones.

\begin{prop}
Let $\B$ be a bialgebra and  $ u $ be its unit. Let $\mathcal{B}'$ be a result of adjoining to $\B$
successive unit elements $e$ and $ 1 $ with respect to the  multiplication and such that:
\begin{eqnarray*}
&&\Delta(1)=1\otimes (e-u)+u\otimes( 1-2e+2u),\\ &&\Delta(e)=e\otimes(e-u)+u\otimes(2u-e), \\
&&\varepsilon(1)=2,\  \varepsilon(e)=2.
\end{eqnarray*} Then $\mathcal{B}'$ is a weak bialgebra.

\end{prop}

\begin{prop}
Let $\B$ be a bialgebra and $u$ its unit. Let $\mathcal{B}'$ be a result of adjoining to $\B$ successive
unit elements $ e $ and $ 1 $ with respect to the  multiplication and  such that:
\begin{eqnarray*}
&&\Delta(1)=(1-e)\otimes (1-e)+(e-u)\otimes(e-u)+u\otimes u,\\ &&\Delta(e)=(e-u)\otimes(e-u)+u\otimes
u,\\ &&\varepsilon(1)=3,\  \varepsilon(e)=2.
\end{eqnarray*}
Then $\mathcal{B}'$ is a weak bialgebra.

\end{prop}

\begin{rem} In the previous propositions, if the bialgebra $\B$ is a Hopf algebra then $\B'$ becomes a weak Hopf algebra by setting $S(1)=1$ and $S(e)=e$.
\end{rem}

\section{Algebraic varieties of weak bialgebras}
Let $V$ be an $n$-dimensional  $\K$-vector space and $\mathfrak{b}=\{e_{1}, \cdots ,e_{n}\}$ be a
basis of $V$.  Let
$\HH=(V,m,\eta,\Delta,\varepsilon)$ (resp. $\HH=(V,m,\eta,\Delta,\varepsilon,S)$) be a  weak bialgebra (resp. weak Hopf algebra). Set  $\eta(1)=e_{1}$ for the unit. Multiplication $m$, comultiplication $\Delta$,  counit $\varepsilon$ and  antipode $S$ write, with respect to this basis,
$$m(e_{i}, e_{j})=\displaystyle\sum_{k=1}^{n} C_{i,j}^{k} e_{k},\ \Delta(e_{k})= \displaystyle\sum_{i,j=1}^{n} D_{k}^{i,j} e_{i}\otimes e_{j}, \ \varepsilon(e_{k})=f_{k},\
 S(e_{i})=\displaystyle\sum_{i,j=1}^{n} s_{i,j}
e_{j}.$$
  The
collection $\{C_{i,j}^{k}, D_{k}^{i,j}, f_{k}\ : \ i, j, k=1, \cdots, n\}$ is the set of structure constants of the weak bialgebra $H$, with
respect to the basis $\mathfrak{b}$. Any $n$-dimensional weak bialgebra is identified to a point of
$\mathbb{K}^{2n^{3}+n}$, determined by a collection $\{C_{i,j}^{k}, D_{k}^{i,j},f_{k} \ : \ i,j,k=1,\cdots,n
\}\in \mathbb{K}^{2n^{3}+n}$, satisfying  the following equations:
\begin{equation}\label{WBeq1}
\displaystyle\sum_{\ell =1}^{n} D_{s}^{\ell ,k} D_{\ell }^{i,j}- D_{s}^{i,\ell }D_{\ell }^{j,k}=0,
\end{equation}
\begin{equation}\label{WBeq2}
\displaystyle\sum_{k=1}^{n}D_{i}^{j,k}f_{k}=\displaystyle
\end{equation}
where $\delta _{i,j}$ is a Kronecker symbol,
\begin{equation}\label{WBeq3}
\displaystyle\sum_{\ell =1}^{n}(C_{i,j}^{\ell }D_{\ell }^{s,r}-\displaystyle \sum_{p,q,t,\ell
=1}^{n}D_{i}^{t,\ell }D_{j}^{p,q}C_{t,p}^{s}C_{\ell ,q}^{r})=0,
\end{equation}
\begin{equation}\label{WBeq4}
\displaystyle\sum_{\ell =1}^{n}(D_{1}^{s,\ell }D_{\ell }^{r,k}-\displaystyle \sum_{p,q,t,\ell
=1}^{n}D_{1}^{p,q}D_{1}^{t,l}C_{1,t}^{s}C_{p,\ell }^{r}C_{q,1}^{k})=0,
\end{equation}
\begin{equation}\label{WBeq5}
\displaystyle\sum_{t,r=1}^{n} C_{i,j}^{t} C_{t,k}^{r}  f_{r}=\displaystyle \sum_{p,q,r,t=1}^{n} D_{j}^{p,q}
C_{q,k}^{r} C_{i,p}^{t} f_{t} f_{r}=\displaystyle \sum_{p,q,r,t=1}^{n} D_{j}^{p,q} C_{i,q}^{r}C_{p,k}^{t}
f_{t} f_{r}.
\end{equation}
We denote by $\mathcal{BF}_{n}$ the set of $n$-dimensional weak bialgebras. The previous system of  equations
endows $\mathcal{BF}_{n}$ with a structure of affine   algebraic variety imbedded in
$\mathbb{K}^{2n^{3}+n}$.

Similarly,
 an  $n$-dimensional weak Hopf algebra $\HH =(V,m,\eta,\Delta,\varepsilon,S)$ is determined, with
respect  to the basis $\mathfrak{b}$ of $V$, by a collection of structure constants
 $\{C_{i,j}^{k}, D_{k}^{i,j},f_{k},s_{i,j}\ : \
i,j,k=1,\cdots,n\}\in\mathbb{K}^{2n^{3}+n^{2}+n}$, satisfying the equations
\eqref{WBeq1}-\eqref{WBeq5}, and moreover the following equations :
\begin{equation}
\displaystyle\sum_{j,r,k=1}^{n}D_{i}^{j,k}s_{k,r}C_{j,r}^{t}-\displaystyle%
\sum_{j,k=1}^{n}D_{1}^{j,t}C_{j,i}^{k}f_{k}=0,
\end{equation}
\begin{equation}
\displaystyle\sum_{j,r,k=1}^{n}D_{i}^{k,j}s_{k,r}C_{r,j}^{t}-\displaystyle%
\sum_{j,k=1}^{n}D_{1}^{t,j}C_{i,j}^{k}f_{k}=0,
\end{equation}
\begin{equation}
\displaystyle\sum_{p,q,j,r,m,\ell ,t=1}^{n}D_{i}^{p,q}D_{p}^{j,r}s_{r,m}s_{q,\ell }C_{m,r}^{t}C_{t,\ell
}^{k}-s_{i,k}=0.
\end{equation}
We denote by $\mathcal{HF}_n $ the set of $n$-dimensional weak Hopf algebras.


We define the action of linear groups on the algebraic varieties of weak bialgebras $\mathcal{BF}_n $ and similarly on the
algebraic varieties of weak  Hopf algebras $\mathcal{HF}_n $:
 $$  GL_{n}(\K) \times
\mathcal{BF}_{n}\rightarrow \mathcal{BF}_{n},$$
$$(g,\HH )\longmapsto g\cdot \HH . $$
This action is defined for all $x,y$ in $V$ by
$$
\begin{array}{l}
(g\cdot m)(x\otimes y)=g^{-1}(m(g(x)\otimes g(y))),\\ (g\cdot\Delta)(x)= g^{-1}\otimes g^{-1}(\Delta
g(x)),\\ (g\cdot\varepsilon)(x)=\varepsilon(g(x)).
\end{array}
$$
The action on the antipode is given by
$$g\cdot S= g^{-1}\circ S\circ g. $$
The orbit of a weak bialgebra (resp. weak Hopf algebra) $\HH$ describes the  isomorphisms class, it is characterized by:
$$\vartheta(\HH)=\{g \cdot \HH :g \in
 GL_{n}(\K)\}.$$ The stabilizer of $\HH $ is  $$  stab(\HH)=\{g \in
 GL_{n}(\K):g\cdot \HH =\HH\},$$ which corresponds to the automorphisms groups of $\HH$.
 We have
 $
 dim \ \vartheta(\HH)=n^{2}-dim \ Aut(\HH).
 $

\section{Classifications and homomorphism groups}
In this section, we establish a classification, up to isomorphism,  of  weak bialgebras
and  weak Hopf algebras of dimension $2$ and $3$.
\subsection{Classification of associative algebras}
The classification of $n$-dimensional associative algebras is known for $n\leq 5$,
(\cite{Gabriel}, \cite{Mazzola}). We recall the results in dimensions $2$ and $3$. Let $\{e_1,\cdots,e_n\}$ be a basis of the underlaying vector space.
\begin{prop}
Every $2$-dimensional associative algebra   is isomorphic to one of the following algebras:
$$m_{1}^{2}(e_{1}, e_{1})=e_{1}, m_{1}^{2}(e_{1}, e_{2})=e_{2}, m_{1}^{2}(e_{2}, e_{1})=e_{2},
m_{1}^{2}(e_{2},
e_{2})=0$$
$$m_{2}^{2}(e_{1}, e_{1})=e_{1}, m_{2}^{2}(e_{1}, e_{2})=e_{2}, m_{2}^{2}(e_{2}, e_{1})=e_{2},
m_{2}^{2}(e_{2},
e_{2})=e_{2}.$$
\end{prop}
\begin{prop} Every  $3$-dimensional associative algebra is isomorphic
to one of the following algebras:\\
$$m_{1}^{3}(e_{1}, e_{1})=e_{1},\ m_{1}^{3}(e_{1}, e_{2})=e_{2}, \ m_{1}^{3}(e_{2}, e_{1})=e_{2}, \
m_{1}^{3}(e_{2}, e_{2})=e_{2},
 \ m_{1}^{3}(e_{1}, e_{3})=e_{3},$$
$$
m_{1}^{3}(e_{3}, e_{1})=e_{3},\ m_{1}^{3}(e_{2}, e_{3})=e_{3},\ m_{1}^{3}(e_{3}, e_{2})=e_{3},\
m_{1}^{3}(e_{3}, e_{3})=e_{3},$$

$$m_{2}^{3}(e_{1},\ e_{1})=e_{1},\ m_{2}^{3}(e_{1}, e_{2})=e_{2},\ m_{2}^{3}(e_{2},
e_{1})=e_{2},\ m_{2}^{3}(e_{2}, e_{2})=e_{2},\ m_{2}^{3}(e_{1}, e_{3})=e_{3},$$
$$
m_{2}^{3}(e_{3}, e_{1})=e_{3},\
 m_{2}^{3}(e_{2}, e_{3})=e_{3},\ m_{2}^{3}(e_{3}, e_{2})=e_{3},\ m_{2}^{3}(e_{3},
 e_{3})=0,$$

$$m_{3}^{3}(e_{1}, e_{1})=e_{1},\ m_{3}^{3}(e_{1}, e_{2})=e_{2},\ m_{3}^{3}(e_{2}, e_{1})=e_{2},
m_{3}^{3}(e_{2}, e_{2})=e_{2},\ m_{3}^{3}(e_{1}, e_{3})=e_{3},$$
$$ m_{3}^{3}(e_{3}, e_{1})=e_{3},\
 m_{3}^{3}(e_{2}, e_{3})=0,\ m_{3}^{3}(e_{3}, e_{2})=0,\ m_{3}^{3}(e_{3},
 e_{3})=0,$$

$$m_{4}^{3}(e_{1}, e_{1})=e_{1},\ m_{4}^{3}(e_{1}, e_{2})=e_{2},\ m_{4}^{3}(e_{2},
e_{1})=e_{2},\ m_{4}^{3}(e_{2}, e_{2})=0,\ m_{4}^{3}(e_{1}, e_{3})=e_{3}, $$
$$m_{4}^{3}(e_{3}, e_{1})=e_{3},\
 m_{4}^{3}(e_{2}, e_{3})=0,\ m_{4}^{3}(e_{3}, e_{2})=0,\ m_{4}^{3}(e_{3},
 e_{3})=0,$$

$$m_{5}^{3}(e_{1}, e_{1})=e_{1},\ m_{5}^{3}(e_{1}, e_{2})=e_{2},\ m_{5}^{3}(e_{2},
e_{1})=e_{2},\ m_{5}^{3}(e_{2}, e_{2})=e_{2},\ m_{5}^{3}(e_{1}, e_{3})=e_{3},$$
$$ m_{5}^{3}(e_{3}, e_{1})=e_{3},\
 m_{5}^{3}(e_{2}, e_{3})=e_{3},\ m_{5}^{3}(e_{3}, e_{2})=0,\ m_{5}^{3}(e_{3},
 e_{3})=0.$$
\end{prop}
Next, we will build  $2$ and $3$-dimensional weak bialgebras and weak Hopf algebras using the previous associative algebras. All the calculations are done using a computer algebra system.

\subsection{Classification of $2$-dimensional  weak bialgebras and weak Hopf algebras}

Let $\{e_{1}, e_{2}\}$ be a basis  of $V=\mathbb{C}^{2}$.

\begin{prop}
 Every $2$-dimensional weak bialgebra is isomorphic to
one of the following weak bialgebras:

\begin{enumerate}\item
$$
\begin{array}{l}m_{2}^{2}(e_{1}, e_{1})=e_{1}, m_{2}^{2}(e_{1}, e_{2})=e_{2}, m_{2}^{2}(e_{2},
e_{1})=e_{2}, m_{2}^{2}(e_{2}, e_{2})=e_{2},\\

\Delta (e_{1} )=e_{1}\otimes e_{1},\\ \Delta (e_{2})=e_{2}\otimes e_{2},\\ \varepsilon (e_1)=1,
\varepsilon (e_2)=1.
\end{array}
$$
\item
$$
\begin{array}{l}
m_{2}^{2}(e_{1}, e_{1})=e_{1}, m_{2}^{2}(e_{1}, e_{2})=e_{2}, m_{2}^{2}(e_{2}, e_{1})=e_{2},
m_{2}^{2}(e_{2}, e_{2})=e_{2},\\ \Delta (e_{1} )=e_{1}\otimes e_{1},\\ \Delta
(e_{2})=(e_{1}-e_{2})\otimes(e_{1}-e_{2})\otimes  +e_{2}\otimes e_{2},\\ \varepsilon (e_1)=1,
\varepsilon (e_2)=1.
\end{array}
$$
\item
$$
\begin{array}{l}
m_{2}^{2}(e_{1}, e_{1})=e_{1}, m_{2}^{2}(e_{1}, e_{2})=e_{2}, m_{2}^{2}(e_{2}, e_{1})=e_{2},
m_{2}^{2}(e_{2}, e_{2})=e_{2},\\ \Delta (e_{1})=(e_{1}-e_{2})\otimes(e_{1}-e_{2})+ e_{2}\otimes
e_{2},\\ \Delta (e_{2})= e_{2}\otimes e_{2},\\ \varepsilon (e_1)=2, \varepsilon (e_2)=1.
\end{array}
$$
\end{enumerate}

\end{prop}

From the precedent classification, we derive $2$-dimensional  weak Hopf algebras.

\begin{prop}

There exist, up to isomorphism, two $2$-dimensional  weak Hopf algebras, which are given by  :
\begin{enumerate}

\item
$$
\begin{array}{l}
m_{2}^{2}(e_{1}, e_{1})=e_{1}, m_{2}^{2}(e_{1}, e_{2})=e_{2}, m_{2}^{2}(e_{2}, e_{1})=e_{2},
m_{2}^{2}(e_{2},e_{2})=e_{2},\\ \Delta (e_{1})= e_{1}\otimes e_{1},  \\ \Delta (e_{2}
)=(e_{1}-e_{2})\otimes(e_{1}-e_{2})+ e_{2}\otimes e_{2},\\ \varepsilon (e_1)=1, \varepsilon
(e_2)=1,\\ S(e_{1})=e_{1}, S(e_{2})=e_{2}.
\end{array}
$$
\item
$$
\begin{array}{l}
m_{2}^{2}(e_{1}, e_{1})=e_{1}, m_{2}^{2}(e_{1}, e_{2})=e_{2}, m_{2}^{2}(e_{2}, e_{1})=e_{2},
m_{2}^{2}(e_{2}, e_{2})=e_{2},\\ \Delta (e_{1} )=( e_{1}-e_{2})\otimes (e_{1}-e_{2})+ e_{2}\otimes
e_{2},\\ \Delta (e_{2})=e_{2}\otimes e_{2},\\ \varepsilon (e_1)=2, \varepsilon (e_2)=1,\\
S(e_{1})=e_{1}, S(e_{2})=e_{2}.
\end{array}
$$
\end{enumerate}
\end{prop}

\subsection{Classification of $3$-dimensional  weak bialgebras and weak Hopf algebras}
  Let $V=\mathbb{C}^{3}$ be a $3$-dimensional vector space with a basis  $\{e_{1}, e_{2},
 e_{3}\}.$ We provide all $3$-dimensional weak bialgebras. Then we specify which of them correspond to  weak Hopf
 algebras.

\begin{prop}Every  $3$-dimensional weak bialgebra  is isomorphic to one
of the following weak bialgebras.
\begin{enumerate}\item

$$\begin{array}{l}
m_{1}^{3}(e_{1}, e_{1})=e_{1}, m_{1}^{3}(e_{1}, e_{2})=e_{2}, m_{1}^{3}(e_{2}, e_{1})=e_{2},
m_{1}^{3}(e_{2}, e_{2})=e_{2}, m_{1}^{3}(e_{1}, e_{3})=e_{3},\\ m_{1}^{3}(e_{3}, e_{1})=e_{3},
m_{1}^{3}(e_{2}, e_{3})=e_{3}, m_{1}^{3}(e_{3}, e_{2})=e_{3}, m_{1}^{3}(e_{3}, e_{3})=e_{3},\\

\Delta (e_{1})=e_{1}\otimes e_{1},\\ \Delta(e_{2})=e_{1}\otimes (e_{1}-e_{3})+ e_{2}\otimes
(2e_{3}-e_{2})+e_{3}\otimes (2e_{2}-e_{3}- e_{1}),\\ \Delta (e_{3} )=e_{1}\otimes
(e_{2}-e_{3})+e_{2}\otimes (e_{1}-2e_{2}+ e_{3})+e_{3}\otimes (e_{2}+e_{3}-e_{1}),
\\ \varepsilon(e _1)=\varepsilon (e_2)=\varepsilon (e_3)=1.
\end{array}
$$

\item
$$
\begin{array}{l}
m_{1}^{3}(e_{1}, e_{1})=e_{1}, m_{1}^{3}(e_{1}, e_{2})=e_{2}, m_{1}^{3}(e_{2}, e_{1})=e_{2},
m_{1}^{3}(e_{2}, e_{2})=e_{2}, m_{1}^{3}(e_{1}, e_{3})=e_{3}, \\ m_{1}^{3}(e_{3}, e_{1})=e_{3},
m_{1}^{3}(e_{2}, e_{3})=e_{3}, m_{1}^{3}(e_{3}, e_{2})=e_{3}, m_{1}^{3}(e_{3}, e_{3})=e_{3},\\

\Delta (e_{1} )=e_{1} \otimes  e_{1},\\ \Delta (e_{2} )=e_{2} \otimes e_{2}, \\ \Delta (e_{3})=e_{3}
\otimes e_{3}, \\ \varepsilon (e_1)=\varepsilon (e_2) = \varepsilon (e_3) = 1.
\end{array}
$$

\item
$$
\begin{array}{l}
m_{1}^{3}(e_{1}, e_{1})=e_{1}, m_{1}^{3}(e_{1}, e_{2})=e_{2}, m_{1}^{3}(e_{2}, e_{1})=e_{2},
m_{1}^{3}(e_{2}, e_{2})=e_{2}, m_{1}^{3}(e_{1}, e_{3})=e_{3},\\ m_{1}^{3}(e_{3}, e_{1})=e_{3},
m_{1}^{3}(e_{2}, e_{3})=e_{3}, m_{1}^{3}(e_{3}, e_{2})=e_{3}, m_{1}^{3}(e_{3}, e_{3})=e_{3},\\

\Delta (e_{1} )=e_{1} \otimes  e_{1}, \\ \Delta (e_{2} )=e_{2} \otimes e_{2}, \\ \Delta
(e_{3})=(e_{2}-e_{3}) \otimes e_{3}+e_{3} \otimes (e_{2}-e_{3}), \\ \varepsilon (e_1) = \varepsilon
(e_2) =1, \varepsilon (e_3) = 0.
\end{array}
$$

\item
$$
\begin{array}{l}
m_{1}^{3}(e_{1}, e_{1})=e_{1}, m_{1}^{3}(e_{1}, e_{2})=e_{2}, m_{1}^{3}(e_{2}, e_{1})=e_{2},
m_{1}^{3}(e_{2}, e_{2})=e_{2}, m_{1}^{3}(e_{1}, e_{3})=e_{3}, \\ m_{1}^{3}(e_{3}, e_{1})=e_{3},
m_{1}^{3}(e_{2}, e_{3})=e_{3}, m_{1}^{3}(e_{3}, e_{2})=e_{3}, m_{1}^{3}(e_{3}, e_{3})=e_{3},\\

\Delta (e_{1} )=e_{1}\otimes  e_{1},\\ \Delta (e_{2} ) =(e_{2}-e_{3}) \otimes e_{3}+e_{3}\otimes
e_{2},\\ \Delta (e_{3} )=e_{3}\otimes e_{3},\\ \varepsilon (e_1) =\varepsilon (e_2) =\varepsilon (e_3)
= 1.
\end{array}
$$

\item
$$
\begin{array}{l}

m_{1}^{3}(e_{1}, e_{1})=e_{1}, m_{1}^{3}(e_{1}, e_{2})=e_{2}, m_{1}^{3}(e_{2}, e_{1})=e_{2},
m_{1}^{3}(e_{2}, e_{2})=e_{2}, m_{1}^{3}(e_{1}, e_{3})=e_{3},\\ m_{1}^{3}(e_{3}, e_{1})=e_{3},
m_{1}^{3}(e_{2}, e_{3})=e_{3}, m_{1}^{3}(e_{3}, e_{2})=e_{3}, m_{1}^{3}(e_{3}, e_{3})=e_{3}, \\

\Delta (e_{1} )=e_{1} \otimes  e_{1},\\ \Delta (e_{2} )=e_{2} \otimes e_{2} + (e_{1}- e_{2} )\otimes
e_{3},\\ \Delta(e_{3})=(e_{1}- e_{3}) \otimes e_{3}+e_{3} \otimes e_{2}, \\ \varepsilon (e_1)
=\varepsilon (e_2)=1,  \varepsilon (e_3)=0.
\end{array}
$$

\item
$$
\begin{array}{l}

m_{1}^{3}(e_{1}, e_{1})=e_{1}, m_{1}^{3}(e_{1}, e_{2})=e_{2}, m_{1}^{3}(e_{2}, e_{1})=e_{2},
m_{1}^{3}(e_{2}, e_{2})=e_{2}, m_{1}^{3}(e_{1}, e_{3})=e_{3},\\ m_{1}^{3}(e_{3}, e_{1})=e_{3},
m_{1}^{3}(e_{2}, e_{3})=e_{3}, m_{1}^{3}(e_{3}, e_{2})=e_{3}, m_{1}^{3}(e_{3}, e_{3})=e_{3}, \\

\Delta (e_{1} )=e_{1} \otimes  e_{1},\\ \Delta (e_{2})=e_{2}\otimes e_{2} + e_{3} \otimes
(e_{1}-e_{2}), \\ \Delta(e_{3})=e_{2} \otimes e_{3} + e_{3} \otimes e_{1}-e_{3} \otimes e_{3},\\
\varepsilon (e_1) = \varepsilon (e_2) =\varepsilon(e_3)= 1.
\end{array}
$$

\item
$$
\begin{array}{l}
m_{1}^{3}(e_{1}, e_{1})=e_{1}, m_{1}^{3}(e_{1}, e_{2})=e_{2}, m_{1}^{3}(e_{2}, e_{1})=e_{2},
m_{1}^{3}(e_{2}, e_{2})=e_{2}, m_{1}^{3}(e_{1}, e_{3})=e_{3},\\ m_{1}^{3}(e_{3}, e_{1})=e_{3},
m_{1}^{3}(e_{2}, e_{3})=e_{3}, m_{1}^{3}(e_{3}, e_{2})=e_{3}, m_{1}^{3}(e_{3}, e_{3})=e_{3},\\

\Delta (e_{1} )=e_{1} \otimes  e_{1}, \\ \Delta (e_{2} )=(e_{1}-e_{2}) \otimes(e_{1} -
e_{2})+e_{2}\otimes e_{2}, \\ \Delta(e_{3})=e_{3} \otimes e_{3},\\ \varepsilon (e_1) = \varepsilon
(e_2) = \varepsilon (e_3) = 1.
\end{array}
$$

\item
$$
\begin{array}{l}
m_{1}^{3}(e_{1}, e_{1})=e_{1}, m_{1}^{3}(e_{1}, e_{2})=e_{2}, m_{1}^{3}(e_{2}, e_{1})=e_{2},
m_{1}^{3}(e_{2}, e_{2})=e_{2}, m_{1}^{3}(e_{1}, e_{3})=e_{3},\\ m_{1}^{3}(e_{3}, e_{1})=e_{3},
 m_{1}^{3}(e_{2}, e_{3})=e_{3}, m_{1}^{3}(e_{3}, e_{2})=e_{3}, m_{1}^{3}(e_{3}, e_{3})=e_{3},\\

\Delta (e_{1} )=(e_{1}-e_{2}) \otimes(e_{1} - e_{2})+e_{2}\otimes e_{2},\\ \Delta (e_{2} )=e_{2}
\otimes e_{2},\\ \Delta (e_{3})=e_{3} \otimes e_{3}, \\ \varepsilon (e_1) =2,  \varepsilon (e_2) = 1,
\varepsilon (e_3) = 1.
\end{array}
$$

\item
$$
\begin{array}{l}
m_{1}^{3}(e_{1}, e_{1})=e_{1}, m_{1}^{3}(e_{1}, e_{2})=e_{2}, m_{1}^{3}(e_{2}, e_{1})=e_{2},
m_{1}^{3}(e_{2}, e_{2})=e_{2}, m_{1}^{3}(e_{1}, e_{3})=e_{3},\\ m_{1}^{3}(e_{3}, e_{1})=e_{3},
 m_{1}^{3}(e_{2}, e_{3})=e_{3}, m_{1}^{3}(e_{3}, e_{2})=e_{3}, m_{1}^{3}(e_{3}, e_{3})=e_{3},\\

\Delta (e_{1} )=(e_{1}-e_{2}) \otimes(e_{1} - e_{2})+e_{2}\otimes e_{2}, \\ \Delta (e_{2} )=e_{2}
\otimes e_{2},\\ \Delta (e_{3})=(e_{2}-e_{3})\otimes e_{3}+e_{3}\otimes(e_{2}-e_{3}), \\ \varepsilon
(e_1) = 2, \varepsilon (e_2) = 1, \varepsilon (e_3) = 0.
\end{array}
$$

\item
$$
\begin{array}{l}
m_{1}^{3}(e_{1}, e_{1})=e_{1}, m_{1}^{3}(e_{1}, e_{2})=e_{2}, m_{1}^{3}(e_{2}, e_{1})=e_{2},
m_{1}^{3}(e_{2}, e_{2})=e_{2}, m_{1}^{3}(e_{1}, e_{3})=e_{3}, \\m_{1}^{3}(e_{3}, e_{1})=e_{3},
 m_{1}^{3}(e_{2}, e_{3})=e_{3}, m_{1}^{3}(e_{3}, e_{2})=e_{3}, m_{1}^{3}(e_{3}, e_{3})=e_{3},\\

\Delta (e_{1} )=(e_{1}-e_{2}) \otimes(e_{1} - e_{2})+(e_{2}-e_{3}) \otimes(e_{2} -
e_{3})+e_{3}\otimes e_{3}, \\ \Delta (e_{2} )=(e_{2}-e_{3}) \otimes(e_{2} - e_{3})+e_{3}\otimes
e_{3}, \\ \Delta (e_{3})=e_{3} \otimes e_{3}, \\ \varepsilon (e_1) = 3, \varepsilon (e_2) =2,
\varepsilon (e_3) = 1.
\end{array}
$$
\item
$$
\begin{array}{l}
m_{1}^{3}(e_{1}, e_{1})=e_{1}, m_{1}^{3}(e_{1}, e_{2})=e_{2}, m_{1}^{3}(e_{2}, e_{1})=e_{2},
m_{1}^{3}(e_{2}, e_{2})=e_{2}, m_{1}^{3}(e_{1}, e_{3})=e_{3}, \\m_{1}^{3}(e_{3}, e_{1})=e_{3},
 m_{1}^{3}(e_{2}, e_{3})=e_{3}, m_{1}^{3}(e_{3}, e_{2})=e_{3}, m_{1}^{3}(e_{3}, e_{3})=e_{3},\\

\Delta(e_{1})=e_{1}\otimes(e_{2}-e_{3})+ e_{3}\otimes(e_{1}-2e_{2}+2e_{3}) , \\ \Delta (e_{2}
)=(e_{2}-e_{3}) \otimes(e_{2} - e_{3})+e_{3}\otimes e_{3},\\ \Delta (e_{3})=e_{3} \otimes e_{3}, \\
\varepsilon (e_1) = 2, \varepsilon (e_2) =2,  \varepsilon (e_3) = 1.
\end{array}
$$

\item
$$
\begin{array}{l}
m_{2}^{3}(e_{1}, e_{1})=e_{1}, m_{2}^{3}(e_{1}, e_{2})=e_{2}, m_{2}^{3}(e_{2}, e_{1})=e_{2},
m_{2}^{3}(e_{2}, e_{2})=e_{2}, m_{2}^{3}(e_{1}, e_{3})=e_{3}, \\m_{2}^{3}(e_{3}, e_{1})=e_{3},
 m_{2}^{3}(e_{2}, e_{3})=e_{3}, m_{2}^{3}(e_{3}, e_{2})=e_{3}, m_{2}^{3}(e_{3}, e_{3})=0,\\

\Delta (e_{1} )=e_{1} \otimes  e_{1}, \\ \Delta (e_{2} )=e_{1} \otimes e_{2} + e_{2} \otimes e_{1}
-e_{2} \otimes e_{2} ,\\ \Delta (e_{3})=e_{1} \otimes e_{3}+e_{3} \otimes e_{1}-e_{3} \otimes e_{2},
\\ \varepsilon (e_1) =1, \varepsilon (e_2) =\varepsilon (e_3) = 0.
\end{array}
$$

\item
$$
\begin{array}{l}
m_{2}^{3}(e_{1}, e_{1})=e_{1}, m_{2}^{3}(e_{1}, e_{2})=e_{2}, m_{2}^{3}(e_{2}, e_{1})=e_{2},
m_{2}^{3}(e_{2}, e_{2})=e_{2}, m_{2}^{3}(e_{1}, e_{3})=e_{3}, \\m_{2}^{3}(e_{3}, e_{1})=e_{3},
 m_{2}^{3}(e_{2}, e_{3})=e_{3}, m_{2}^{3}(e_{3}, e_{2})=e_{3}, m_{2}^{3}(e_{3}, e_{3})=0,\\

\Delta (e_{1} )=e_{1} \otimes  e_{1}, \\ \Delta (e_{2} )=e_{1} \otimes e_{2} + e_{2} \otimes e_{1}
-e_{2} \otimes e_{2} ,\\ \Delta (e_{3})=e_{1} \otimes e_{3}-e_{2} \otimes e_{3}+e_{3} \otimes
e_{1}-e_{3} \otimes e_{2}+e_{3} \otimes e_{3}, \\ \varepsilon (e_1) = 1, \varepsilon (e_2)
=\varepsilon (e_3) = 0.
\end{array}
$$

\item
$$
\begin{array}{l}
m_{2}^{3}(e_{1}, e_{1})=e_{1}, m_{2}^{3}(e_{1}, e_{2})=e_{2}, m_{2}^{3}(e_{2}, e_{1})=e_{2},
m_{2}^{3}(e_{2}, e_{2})=e_{2}, m_{2}^{3}(e_{1}, e_{3})=e_{3}, \\m_{2}^{3}(e_{3}, e_{1})=e_{3},
 m_{2}^{3}(e_{2}, e_{3})=e_{3}, m_{2}^{3}(e_{3}, e_{2})=e_{3}, m_{2}^{3}(e_{3}, e_{3})=0,\\

\Delta (e_{1} )=e_{1} \otimes  e_{1}, \\ \Delta (e_{2} )=e_{1} \otimes e_{2} + e_{2} \otimes e_{1}
-e_{2} \otimes e_{2} ,\\ \Delta (e_{3})=e_{1} \otimes e_{3}-e_{2} \otimes e_{3}+e_{3} \otimes e_{1},
\\ \varepsilon (e_1) = 1, \varepsilon (e_2) =\varepsilon (e_3) = 0.
\end{array}
$$

\item
$$
\begin{array}{l}
m_{2}^{3}(e_{1}, e_{1})=e_{1}, m_{2}^{3}(e_{1}, e_{2})=e_{2}, m_{2}^{3}(e_{2}, e_{1})=e_{2},
m_{2}^{3}(e_{2}, e_{2})=e_{2}, m_{2}^{3}(e_{1}, e_{3})=e_{3}, \\ m_{2}^{3}(e_{3}, e_{1})=e_{3},
 m_{2}^{3}(e_{2}, e_{3})=e_{3}, m_{2}^{3}(e_{3}, e_{2})=e_{3}, m_{2}^{3}(e_{3}, e_{3})=0,\\

\Delta (e_{1} )=e_{1} \otimes  e_{1}, \\ \Delta (e_{2} )=e_{1} \otimes e_{2} + e_{2} \otimes e_{1}
-e_{2} \otimes e_{2} ,\\ \Delta (e_{3})=e_{1} \otimes e_{3}-e_{3} \otimes e_{2}+e_{3} \otimes
e_{1}-e_{2} \otimes e_{3}, \\ \varepsilon (e_1) = 1, \varepsilon (e_2) =\varepsilon (e_3) = 0.
\end{array}
$$

\item
$$
\begin{array}{l}
m_{3}^{3}(e_{1}, e_{1})=e_{1}, m_{3}^{3}(e_{1}, e_{2})=e_{2}, m_{3}^{3}(e_{2}, e_{1})=e_{2},
m_{3}^{3}(e_{2}, e_{2})=e_{2}, m_{3}^{3}(e_{1}, e_{3})=e_{3}, \\m_{3}^{3}(e_{3}, e_{1})=e_{3},
 m_{3}^{3}(e_{2}, e_{3})=e_{3}, m_{3}^{3}(e_{3}, e_{2})=0, m_{3}^{3}(e_{3}, e_{3})=0,\\

\Delta (e_{1} )=(e_{1}-e_{2}) \otimes  (e_{1}-e_{2})+e_{2} \otimes e_{2}, \\ \Delta (e_{2} )=e_{2}
\otimes e_{2}, \\ \Delta (e_{3})=e_{3} \otimes e_{3}, \\ \varepsilon (e_1) =2, \varepsilon (e_2)
=\varepsilon (e_3) = 1.
\end{array}
$$

\item
$$
\begin{array}{l}
m_{3}^{3}(e_{1}, e_{1})=e_{1}, m_{3}^{3}(e_{1}, e_{2})=e_{2}, m_{3}^{3}(e_{2}, e_{1})=e_{2},
m_{3}^{3}(e_{2}, e_{2})=e_{2}, m_{3}^{3}(e_{1}, e_{3})=e_{3}, \\m_{3}^{3}(e_{3}, e_{1})=e_{3},
 m_{3}^{3}(e_{2}, e_{3})=e_{3}, m_{3}^{3}(e_{3}, e_{2})=0, m_{3}^{3}(e_{3}, e_{3})=0,\\

\Delta (e_{1} )=e_{1} \otimes  e_{1}, \\ \Delta (e_{2} )=e_{1} \otimes e_{2}+e_{2} \otimes e_{1}-e_{2}
\otimes e_{2}, \\ \Delta (e_{3})=e_{1} \otimes e_{3}+e_{3} \otimes e_{1}-e_{2} \otimes e_{3}-e_{3}
\otimes e_{2}, \\ \varepsilon (e_1) =1, \varepsilon (e_2) =\varepsilon (e_3) = 0.
\end{array}
$$

\item
$$
\begin{array}{l}
m_{3}^{3}(e_{1}, e_{1})=e_{1}, m_{3}^{3}(e_{1}, e_{2})=e_{2}, m_{3}^{3}(e_{2}, e_{1})=e_{2},
m_{3}^{3}(e_{2}, e_{2})=e_{2}, m_{3}^{3}(e_{1}, e_{3})=e_{3}, \\m_{3}^{3}(e_{3}, e_{1})=e_{3},
 m_{3}^{3}(e_{2}, e_{3})= e_{3}, m_{3}^{3}(e_{3}, e_{2})=0, m_{3}^{3}(e_{3}, e_{3})=0,\\

\Delta (e_{1} )=e_{1} \otimes  e_{1}, \\ \Delta (e_{2} )=e_{1} \otimes e_{2}+e_{2} \otimes e_{1}-e_{2}
\otimes e_{2}-e_{3} \otimes e_{3}, \\ \Delta (e_{3})=e_{1} \otimes e_{3}-e_{2} \otimes e_{3}+e_{3}
\otimes e_{1}-e_{3} \otimes e_{2}, \\ \varepsilon (e_1) =1, \varepsilon (e_2) =\varepsilon (e_3) = 0.
\end{array}
$$

\item
$$
\begin{array}{l}
m_{3}^{3}(e_{1}, e_{1})=e_{1}, m_{3}^{3}(e_{1}, e_{2})=e_{2}, m_{3}^{3}(e_{2}, e_{1})=e_{2},
m_{3}^{3}(e_{2}, e_{2})=e_{2}, m_{3}^{3}(e_{1}, e_{3})=e_{3}, \\m_{3}^{3}(e_{3}, e_{1})=e_{3},
 m_{3}^{3}(e_{2}, e_{3})=e_{3}, m_{3}^{3}(e_{3}, e_{2})=0, m_{3}^{3}(e_{3}, e_{3})=0,\\

\Delta (e_{1} )=e_{1} \otimes  e_{1}, \\ \Delta (e_{2} )=e_{2} \otimes e_{2}, \\ \Delta (e_{3})=e_{2}
\otimes e_{3}+e_{3} \otimes e_{2}, \\ \varepsilon (e_1) =\varepsilon (e_2) =1, \varepsilon (e_3) =
0.

\end{array}
$$

\item
$$
\begin{array}{l}
m_{5}^{3}(e_{1}, e_{1})=e_{1}, m_{5}^{3}(e_{1}, e_{2})=e_{2}, m_{5}^{3}(e_{2}, e_{1})=e_{2},
m_{5}^{3}(e_{2}, e_{2})=e_{2}, m_{5}^{3}(e_{1}, e_{3})=e_{3},\\ m_{5}^{3}(e_{3}, e_{1})=e_{3},
 m_{5}^{3}(e_{2}, e_{3})=e_{3}, m_{5}^{3}(e_{3}, e_{2})=0, m_{5}^{3}(e_{3}, e_{3})=0,\\

\Delta (e_{1} )=e_{1} \otimes  e_{1}, \\ \Delta (e_{2} )=e_{2} \otimes e_{2}+e_{3} \otimes e_{3}, \\
\Delta (e_{3})=e_{2} \otimes e_{3}+e_{3} \otimes e_{2}, \\ \varepsilon (e_1) =\varepsilon (e_2) =1,
\varepsilon (e_3) = 0.

\end{array}
$$

\end{enumerate}

\end{prop}
The $3$-dimensional weak Hopf algebras  are given by the following proposition:

\begin{prop}
Every $3$-dimensional weak Hopf algebra   is isomorphic to one of the  following weak Hopf
algebras.
\begin{enumerate}

\item
$$
\begin{array}{l}
m_{1}^{3}(e_{1}, e_{1})=e_{1}, m_{1}^{3}(e_{1}, e_{2})=e_{2}, m_{1}^{3}(e_{2}, e_{1})=e_{2},
m_{1}^{3}(e_{2}, e_{2})=e_{2}, m_{1}^{3}(e_{1}, e_{3})=e_{3},\\ m_{1}^{3}(e_{3}, e_{1})=e_{3},
 m_{1}^{3}(e_{2}, e_{3})=e_{3}, m_{1}^{3}(e_{3}, e_{2})=e_{3}, m_{1}^{3}(e_{3}, e_{3})=e_{3},\\

\Delta (e_{1} )=e_{1} \otimes  e_{1}, \\ \Delta (e_{2} )=e_{1} \otimes e_{2} + e_{2} \otimes e_{1} -
e_{2} \otimes e_{2} - e_{2} \otimes e_{3} - e_{3} \otimes e_{2} + 2e_{3} \otimes e_{3},\\ \Delta
(e_{3})=e_{1} \otimes e_{3} + e_{2} \otimes e_{2} - 2e_{2} \otimes e_{3} + e_{3} \otimes e_{1} -
2e_{3} \otimes e_{2} + e_{3} \otimes e_{3},\\ \varepsilon (e_1) = 1, \varepsilon (e_2) = \varepsilon
(e_3) = 0,\\ S(e_{1})=e_{1}, S(e_{2})=e_{2}, S(e_{3})=e_{2} - e_{3}.
\end{array}
$$
\item
$$
\begin{array}{l}
m_{1}^{3}(e_{1}, e_{1})=e_{1}, m_{1}^{3}(e_{1}, e_{2})=e_{2}, m_{1}^{3}(e_{2}, e_{1})=e_{2},
m_{1}^{3}(e_{2}, e_{2})=e_{2}, m_{1}^{3}(e_{1}, e_{3})=e_{3},\\ m_{1}^{3}(e_{3}, e_{1})=e_{3},
 m_{1}^{3}(e_{2}, e_{3})=e_{3}, m_{1}^{3}(e_{3}, e_{2})=e_{3}, m_{1}^{3}(e_{3}, e_{3})=e_{3},\\

\Delta (e_{1} )=(e_{1}-e_{2})\otimes (e_{1} - e_{2})+e_{2}\otimes e_{2},\\ \Delta (e_{2}
)=e_{2}\otimes e_{2}, \\ \Delta(e_{3})=(e_{2}-e_{3})\otimes e_{3}+e_{3} \otimes (e_{2}-e_{3}),\\
\varepsilon (e_1) = 2, \varepsilon (e_2) =1, \varepsilon (e_3) = 0,\\ S(e_{1})=e_{1}, S(e_{2})=e_{2},
S(e_{3})=e_{3}.
\end{array}
$$

\item
$$
\begin{array}{l}
m_{1}^{3}(e_{1}, e_{1})=e_{1}, m_{1}^{3}(e_{1}, e_{2})=e_{2}, m_{1}^{3}(e_{2}, e_{1})=e_{2},
m_{1}^{3}(e_{2}, e_{2})=e_{2}, m_{1}^{3}(e_{1}, e_{3})=e_{3},\\ m_{1}^{3}(e_{3}, e_{1})=e_{3},
 m_{1}^{3}(e_{2}, e_{3})=e_{3}, m_{1}^{3}(e_{3}, e_{2})=e_{3}, m_{1}^{3}(e_{3}, e_{3})=e_{3},\\

\Delta (e_{1} )=(e_{1}-e_{2})\otimes (e_{1} - e_{2})+(e_{2}-e_{3})\otimes (e_{2} -
e_{3})+e_{3}\otimes e_{3}, \\ \Delta (e_{2} )=(e_{2}-e_{3})\otimes (e_{2} - e_{3})+e_{3}\otimes
e_{3},\\ \Delta (e_{3})= e_{3} \otimes e_{3},
\\ \varepsilon (e_1) = 3, \varepsilon (e_2) =2,  \varepsilon (e_3) = 1,\\
S(e_{1})=e_{1}, S(e_{2})=e_{2}, S(e_{3})= e_{3}.
\end{array}
$$

\end{enumerate}
\end{prop}

\subsection{Automorphisms group}
In this Section, we compute the automorphisms groups of $2$-dimensional and $3$-dimensional  weak
bialgebras  and weak Hopf algebras obtained above. First, we write down the
conditions which should be satisfied in order that two weak bialgebras lie in the same orbit.

 Let $\HH_{1}=(V, m, \Delta_{1}, \varepsilon_{1})$ and
$ \HH_{2}=(V, m, \Delta_{2}, \varepsilon_{2})$ be two weak bialgebras
 with the same orbit, then there exists  linear bijective map
 $g:V\rightarrow V$ which ensure the transport of the structure. We set,
 with respect to a basis $\{e_{i}\}_{i=1,\cdots,n}$,
 \begin{eqnarray*}
 && g(e_{i})=\displaystyle \sum_{j =1}^{n} T_{i,j} e_{j},\quad
 m(e_{i},e_{j})=\displaystyle\sum_{k=1}^{n} C_{i,j}^{k} e_{k},\\&&
 \Delta_{1}(e_{i})= \displaystyle \sum_{j,k =1}^{n}D_{1,i}^{j,k} e_{j}\otimes e_{k}, \quad\Delta_{2}(e_{i})= \displaystyle \sum_{j,k =1}^{n}D_{2,i}^{j,k} e_{j}\otimes e_{k}, \\&& \varepsilon_{1}(e_{i})= f_{1,i},\quad \varepsilon_{2}(e_{i})=
 f_{2,i}.
 \end{eqnarray*}
 These two weak bialgebras $\HH_1$ and $\HH_2$ are isomorphic if the following conditions are satisfied:

\begin{equation} \displaystyle \sum_{p =1}^{n} T_{i,p}
D_{1,p}^{s,r}-\displaystyle \sum_{p,q =1}^{n} D_{2,i}^{p,q} T_{p,s} T_{q,r}=0     \quad   i,s,r=1,\cdots,n,
\end{equation}
\begin{equation}\displaystyle\sum_{j =1}^{n} T_{i,j}f_{1,j}-f_{2,i}=0 \quad i=1,\cdots,n,
\end{equation}

\begin{equation}\displaystyle\sum_{t =1}^{n} T_{t,k} C_{i,j}^{t}-\displaystyle
\sum_{s,r =1}^{n} T_{i,s} T_{j,r} C_{s,r}^{k}=0 \quad i,j,k=1,\cdots,n.
\end{equation}
\subsubsection{Automorphisms group of 2-dimensional weak bialgebras}
The automorphisms groups of all 2-dimensional weak bialgebras are groups of order $2$
given by :

$$ G=<\{\left(
\begin{array}{cc}
1  &   1 \\ 0  &   -1 \\
\end{array}
\right)\}> .
$$
We obtain a similar group for 2-dimensional weak Hopf algebras.
\subsubsection{Automorphisms group of 3-dimensional weak bialgebras}

The automorphisms group of the weak bialgebras  
$(1), (2), (3), (4), (5), (6), (7), (8),(9), (10), (11), $ is the group of order $6$ given by: $$G=<\{\left(
\begin{array}{ccc}
1  &   0  &   0 \\ 0  &   1  &   1 \\ 0  &   0  &   -1 \\
\end{array}
\right),\left(
\begin{array}{ccc}
1  &   1  &   0 \\ 0  &   0  &   1 \\ 0  &   -1  &   -1 \\
\end{array}
\right)\}> .$$

The automorphisms group of the weak bialgebras (12), (13), (14), (15), (16), (17), (19), is the group :
$$G=\{\left(
\begin{array}{ccc}
1 & 0 & 0\\ 0 & 1 & 0\\ 0 & 0 & \alpha^{\theta}  \\
\end{array}
\right), \theta\in\mathbb{Z},  \alpha \in\mathbb{C}^{\ast}\}. $$

The automorphisms group of the weak bialgebra (18) is the group : $$G=<\{\left(
\begin{array}{ccc}
1 &  0   & 0\\ 0 &  1   & 0\\ 0 & -r/2 & \pm1/2\sqrt[2]{(4e-r^{2})}\\
\end{array}
\right), 4e-r^{2}\neq0, r, e\in\mathbb{C}\}>. $$

The automorphisms group of the weak bialgebra (20) is the group :
$$G=<\{\left(
\begin{array}{ccc}
1 &  0   & 0\\ 0 &  1   & 0\\ 0 & r/2  & \pm1/2\sqrt[2]{(4e+r^{2})}\\
\end{array}
\right), 4e+r^{2}\neq0, r, e\in\mathbb{C}\}>. $$ The automorphisms group of the $3$-dimensional weak Hopf algebras   is the group of order $6$ given by  : $$G=<\{\left(
\begin{array}{ccc}
1  &   0  &   0 \\ 0  &   1  &   1 \\ 0  &   0  &   -1 \\
\end{array}
\right),\left(
\begin{array}{ccc}
1  &   1  &   0 \\ 0  &   0  &   1 \\ 0  &   -1  &   -1 \\
\end{array}
\right)\}> .$$

\end{document}